\documentclass[a4,12pt]{amsart}
\usepackage{amsmath}
\usepackage{amssymb}
\begin{document}
\title{Interior regularity on the Abreu equation}
\author{ Bohui Chen$^1$}
\address{ Yangtze Center of Mathematics, Department of Mathematics, Sichuan University\\
        Chengdu}
\email{bohui@cs.wisc.edu}

 \author  { An-Min Li$^{2}$}
\address{
 Yangtze Center of Mathematics, Department of Mathematics, Sichuan University\\
        Chengdu}
\email{math$\_$li@yahoo.com.cn}

\author{Li Sheng$^{3}$}
\address{Dept. of Math. Sichuan University}
\email{l\_sheng@yahoo.cn }
\footnote[1]{partially supported by NKBRPC(2006CB805905), NSFC
10631050 }
 \footnote[2]{partially supported by
NKBRPC(2006CB805905), NSFC 10631050  and RFDP(20060610004)}
\footnote[3]{Corresponding author}
\renewcommand{\thefootnote}{\fnsymbol{footnote}}

%\headheight 0in

%\headsep 0in

%\footskip 0.5in

%\footheight 0in

\numberwithin{equation}{section}

\newtheorem{theorem}{Theorem}[section]
\newtheorem{assertion}[theorem]{Assertion}
\newtheorem{claim}[theorem]{Claim}
\newtheorem{conjecture}[theorem]{Conjecture}
\newtheorem{corollary}[theorem]{Corollary}
\newtheorem{defn}[theorem]{Definition}
\newtheorem{example}[theorem]{Example}
\newtheorem{figger}[theorem]{Figure}
\newtheorem{lemma}[theorem]{Lemma}
\newtheorem{prop}[theorem]{Proposition}
\newtheorem{remark}[theorem]{Remark}
\newtheorem{assumption}[theorem]{Assumption}

 \def \J{{\cal J}}
              \def \Map{Map(S^2, V)}
              \def \M{{\mathcal M}}
              \def \A{{\mathcal A}}
              \def \B{{\mathcal B}}
              \def \C{{\bf C}}
              \def \Z{{\bf Z}}
              \def \R{{\bf R}}
              \def \P{{\bf P}}
              \def \I{{\bf I}}
              \def \N{{\bf N}}
              \def \T{{\bf T}}
              \def \O{{\mathcal O}}
              \def \Q{{\bf Q}}
              \def \D{{\bf D}}
              \def \H{{\bf H}}
              \def \s{{\mathcal S}}
              \def \e{{\bf E}}
              \def \k{{\bf k}}
              \def \U{{\mathcal U}}
              \def \E{{\mathcal E}}
              \def \F{{\mathcal F}}
              \def \L{{\mathcal L}}
              \def \K{{\mathcal K}}
              \def \G{{\bf G}}

\def \mc{\mathcal}
\def \v{\vskip 0.1in}
\def \n{\noindent}

\def \mfk{\mathfrak}

\def \inv{^{-1}}

\def \real{\mathbb{R}}

\def \upalpha{^{(\alpha)}}
\def \funcR{\mathfrak{R}}

\maketitle

\abstract In this paper we prove the interior regularity for the
solution to the Abreu equation in any dimension assuming the
existence of the $C^0$ estimate.
\endabstract

\vskip 0.1in \noindent \vskip 0.1in \noindent

\section{Introduction}\label{sect_1}
\vskip 0.1in \noindent

Studying extremal metrics within a given K\"ahler class on K\"ahler
manifolds is one of the most important problems in complex
geometry.  Since Calabi introduced the notion of extremal metrics, it  has been studied
intensively in the past decades. There are three
aspects of the problem: sufficient conditions, necessary
conditions and the uniqueness of the existence.  The necessary
conditions for the existence are conjectured to be related to
certain stabilities. There are many works on this aspect
(\cite{T}, \cite{D1},\cite{D2},\cite{CT}). For example, in
\cite{D1}, Donaldson proved that the existence of K\"ahler metrics
of constant scalar curvature implies the Chow or Hilbert
(semi)stability. The uniqueness  is completed by Donaldson, Mabuchi and
Chen-Tian. In \cite{D1}, Donaldson proved that K\"ahler metrics of
constant scalar curvature are unique in any rational K\"ahler class
on any projective manifolds without nontrivial holomorphic vector
fields. The complete answer is answered by Mabuchi for the algebraic case and
answered  by Chen-Tian(\cite{CT}) for general cases. They
showed  that K\"ahler metrics of constant scalar curvature are unique
in any rational K\"ahler class on any compact K\"ahler manifold.

On the other hand, there has been few progresses on the existence
of extremal metrics or K\"ahler metrics of constant scalar curvature.
One reason is that the equation is highly nonlinear and of 4th
order. The general expectation of the problem may be
stated in the
following conjecture ( see \cite{D3}): \vskip 0,1in \noindent {\bf
Conjecture} A smooth polarised projective variety $(V,L)$ admits a
K\"ahler metric of constant scalar curvature  in the class $c_1(L)$ if
and only if it is $K$-stable. \vskip 0,1in \noindent

Donaldson has initiated a program on the extremal metrics on toric varieties.
Since a $2n$-dimensional toric variety can be represented by a convex polytope
in $\real^n$, the problem can be formulated on this polytope.
   In \cite{Ab1},  using Gullimin's method(\cite{G}), Abreu reduced the equation for
K\"ahler metrics of constant scalar curvature
  to an equation on  the polytope.
 The  equation   is now called
the Abreu equation.
On the other hand, Donaldson(\cite{D3})  formulated the problem as an variational problem
on the polytope and related it to the $K$-stability.
Hence, the problem is to solve the Abreu equation with the assumption of the stability condition.
Since the equation is a degenerate 4th order  PDE,  it is a very hard problem and  the
progress  is slow.
We list some important progresses obtained by Donaldson: (1) in \cite{D3} Donaldson
introduces a stronger version of stability (cf. Proposition 5.1.2 and Proposition 5.2.2 in \cite{D3});
Donaldson proved the
existence of weak solutions when this stronger stability holds;
(2) in \cite{D4}, he proved the interior estimates for the Abreu equations on toric surfaces;
(3) in \cite{D5}, he proved the conjecture for the toric surfaces under the assumption of $M$-condition,
and (4) in \cite{D6}, he completely solved the conjecture for the
constant  scalar curvature metrics on toric surfaces.
On the other hand, in \cite{ZZ},Zhou-Zhu proved the existence of weak solutions under certain
properness condition.

In this paper, we prove the interior estimates for the solution
to the Abreu equation in any dimension assuming the existence of the $C^0$ estimate.
Be precisely, the statement is
\begin{theorem}\label{main_thm}
Let $M$ be a real $2n$-dimensional compact toric manifold and $\Delta$ be its Delzant
polytope. Let $K\in C^{\infty}(\overline{\Delta})$. If $u$ is a
solution to the Abreu equation $\s(u)=K$ (cf. \eqref{eqn_1.1}), i.e,
it gives a metric on toric variety of curvature $K$, then for any
domain $\Omega\subset \subset \Delta$ and positive integer $k$,
there exists a constant $C$ depending on $n,|K|_{C^{m}(\overline
\Delta)}, |u|_{C^0(\overline \Delta)}$ and $d:=dist(\Omega,\partial \Delta)$ such
that
\begin{equation}\label{main_est}
|u|_{C^{k,\alpha}(\Omega)}\leq C(n, |K|_{C^m(\overline{\Delta})},
|u|_{C^0(\overline \Delta)}, d),
\end{equation}
where   $m=\max\{n-3,0\}.$
\end{theorem}

{\bf Acknowledgements.} The authors would like to thank Prof.
Yongbin Ruan, Prof.Xiuxiong Chen for many valuable discussions.

\section{Toric manifolds and the Abreu equation}\label{sect_2}

A toric manifold is an $2n$-dimensional symplectic manifold $(M,\omega)$
admitting a $T^n$ Hamiltonian action. The image of the moment map
of the Hamiltonian action is then a polytope in $\real ^n$. Such polytopes
are called Delzant polytopes.
Be precisely, we have the following
definition.
\def \integer{\mathbb{Z}}
\begin{defn}\label{defn_1.1}
A convex polytope $\Delta$ in $\mathbb{R}^n$ is a Delzant polytope
if
\begin{enumerate}
\item there are n edges meeting at each vertex ;
 \item the edges meeting at the vertex $p$ are
rational, i.e. each edge is of the form $p+tv_i$, $0\leq t\leq
\infty$, where $v_i \in \integer^n$;
 \item $v_1,...,v_n$ described above can be chosen
to form a basis of $\integer^n$.
 \end{enumerate}
 \end{defn}
A Delzant polytope can be  defined by $d$ linear inequalities
$$
\langle\xi, h_k\rangle \geq \lambda_k,\;\;\; k = 1,\cdots,d,
$$
where $h_k's$ are primitive elements of the lattice $\mathbb{Z}^n$.
Set $\ell_i: \mathbb{R}^n\rightarrow \real$ be the function $\ell_i
(\xi) = \langle\xi,h_i\rangle - \lambda_i,$ and let $\Delta^{o}$ be
the interior of $\Delta$. In this paper, we assume that
$\Delta={\Delta^o}$. Then $\xi \in \Delta$ if and only if
$\ell_i(\xi)>0$ for all i.

In [De] Delzant associates
to every Delzant polytope $\Delta$ a closed connected symplectic toric
manifold $(M_{\Delta}, \omega)$ of dimension $2n$ together with a
Hamiltonian $T^n$-action
 whose moment map $\varphi: M_\Delta\rightarrow R^n$
satisfies
 $\varphi (M_\Delta)=\overline{\Delta}$. On $(M_{\Delta},
\omega)$, Guillemin shows that  there is
 a natural K\"ahler metric called the Guillemin metric,
denoted by $G_g$ or $\omega_g$ (cf. \eqref{eqn_gui}).

Let $T=T^n$ and $C^\infty_T(M_\Delta)$ be the set of $T$-invariant functions.
Denote  \def \cplane{\mathbb{C}}
$$\M = \{\phi \in C^{\infty}_T(M_\Delta)|\omega_\phi = \omega_g +
\frac{\sqrt{-1}}{2\pi}\partial\bar{\partial}\phi >0\}.$$ For any
$\phi \in\M$, restricting on $\varphi^{-1}(\Delta)$, we have a
K\"ahler potential
$$f(x)=g(x) + \phi(x),$$
where $g(x)$ is the Guillemin potential function for $\omega_g$.
There are natural isomorphisms \def \inv{^{-1}}
$$
\varphi\inv(\Delta)\cong (\cplane^\ast)^n\cong (\real\times
S^1)^n\cong \real^n\times T,
$$
then $f$ and $g$ are treated as functions on $\real^n$ whose coordinates
are still denoted by $x_1,\ldots,x_n$. By the Legendre transformation, we
 put
$$\xi_i=\frac{\partial f}{\partial x_i},\;\;\;\;u(\xi)= \sum \xi_i\frac{\partial f}{\partial x_i} - f(x).$$
It is known that $u(\xi)$ is a smooth function on $\Delta$.
Similarly, the Legendre transformation of $g$ also defines a smooth
function $v$ on $\Delta$. By Guillemin's theorem, $v$ is given by
\begin{equation}\label{eqn_gui}
v(\xi)=\sum_{k=1}^d\ell_k(\xi)\log\ell_k(\xi) .
\end{equation}
For each $u$, set $\psi_u=u-v$. It is known that $\psi_u$ is a
smooth function on $\overline{\Delta}$.

We will use $\s(f)$
or $\s(u)$, or even $\s(\phi)$ to denote the scalar
curvature for the
K\"ahler metric $ \omega_\phi $. The Abreu equation(see \cite{D1}) is
\begin{equation}\label{eqn_1.1}
\s(u)=\sum_{i,j=1}^n U^{ij} w_{ij} = - K\;\;\;in\;\Delta
\end{equation}
where $K$ is a given function on $\overline{\Delta } $, $U^{ij}$ is
the cofactor of the Hessian matrix $D^2 u$ of the convex function
$u$ and
\begin{equation*}\label{eqn_1.2}
w =[\det D^2 u]^{-1}.
\end{equation*}
Set
$$
\funcR(\Delta)=\{u=v+\psi|\psi\in \C^\infty(\overline{\Delta})\}
$$
and
$$
\funcR(\Delta, b)=\{u\in \funcR(\Delta)||\s(u)|\leq b\}.
$$
Let $u$ be a function on $\Delta$. We say it is normal at $p\in
\Delta$ if $u(p)=0$ and $\nabla(p)=0$. It is known that for any
function and any point $p\in \Delta$, the function can be normalized
by adding an affine function. Set $\funcR_p(\Delta)\subset
\funcR(\Delta)$,
 ($\funcR_p(\Delta,b)\subset \funcR(\Delta,b)$ )
 be the set of normal functions at $p$.
\begin{remark}
In this paper, $K$ can be any smooth function.
We remark that
the K\"ahler metric is extremal  if and only if $K$
is an affine function in $\xi$.
\end{remark}

\section { Determinant estimates }\label{sect_3}
The key estimate in this paper is Proposition \ref{prop_3.1}, which
holds for any K\"ahler manifold other than toric manifolds.

  Let $(M,G)$ be a compact complex
manifold with a K\"ahler metric $G$. Let $\omega_G$ be its K\"ahler
form. Denote
$$\widetilde{\M} = \{\phi \in C^{\infty}(M)| \omega_\phi = \omega_G +
\frac{\sqrt{-1}}{2\pi}\partial\bar{\partial}\phi >0\}.$$ Choose a
local coordinate system $z_1,...,z_n$, let $g$ be a local potential
function of $G$ in this coordinate system. For any $\phi\in \widetilde\M$, let
$f = g + \phi$. Denote
$$
H:=\frac{\det(g_{i\bar{j}})}{\det(f_{i\bar{j}})}
,$$
it is known that $H$ is a global function defined on $M$.
 Set
%\begin{equation}
%Tr(f,g) = \sum f^{i\bar{j}}g_{i\bar{j}}, \end{equation}
$
\K= \max\limits_{M}\|Ric(g)\| _{g},
$ where $Ric(g)$ denotes Ricci tensor for metric $G$.
We prove
\begin{prop}\label{prop_3.1}  For any $\phi\in \widetilde\M$ we have
\begin{equation}\label{eqn_3.3}
H \leq \left(2+\frac{\max_M|\s(\phi)|}{n(\K+1)}\right)^n\exp
\left\{(2\K+1)(\max_M\{\phi\}-\min_M\{\phi\})\right\}.
\end{equation}
\end{prop}
 {\bf Proof.} Consider the function
$$\F := \exp\{-C\phi\}H,$$
where $C$ is a constant to be determined later. $\F$ attains its
maximum at a point $p^*\in M$. We have, at $p^*$,
$$- C f^{i\bar j} \phi_{i\bar j}+f^{i\bar j}(\log H)_{i\bar j}\leq 0
$$
which implies
\begin{equation}\label{eqn_3.2}
- C f^{i\bar j} \phi_{i\bar j} +  \s (\phi)+
\sum f^{i\bar{j}}\left(\log
\det(g_{k\bar{l}})\right)_{i\bar{j}} \leq 0.
\end{equation}
We can choose the
coordinates $z_1,...,z_n$ at a neighborhood of $p^\ast$ such that
$$f_{i\bar{j}}= \lambda_i \delta_{ij},\;\;\;g_{i\bar{j}}= \mu_i
\delta_{ij}.$$ From \eqref{eqn_3.2}, we get
$$C\left( \tfrac{\mu_1}{\lambda_1} + \cdot\cdot\cdot +
\tfrac{\mu_n}{\lambda_n}\right) - Cn +  \s - \left(
\tfrac{\mu_1}{\lambda_1} + \cdot\cdot\cdot +
\tfrac{\mu_n}{\lambda_n}\right) \K \leq 0.$$ We choose $C = 2\K+1$ and
apply an elementary inequality
$$
\frac{1}{n}\left(\tfrac{\mu_1}{\lambda_1} + \cdot\cdot\cdot +
\tfrac{\mu_n}{\lambda_n}\right)\geq \left(\frac{\mu_1\cdots\mu_n}{\lambda_1\cdots
\lambda_n}\right)^{1/n}
$$
to get
$$n(\K+1) H^{\tfrac{1}{n}} \leq n(2\K+1) +  |\s|.$$
It follows that, at $p^*$,
$$\exp\{-C\phi\}H \leq \left(2+\frac{|\s|}{n(\K+1)}\right)^n\exp\{-(2\K+1)\min_M\{\phi\}\}.$$
Then \eqref{eqn_3.3} follows. $\Box$

\vskip 0.1in On the other hand, on $\Delta$, Donaldson gives
several estimates for $u\in \funcR(\Delta, b)$.

\begin{lemma}\label{lemma_3.2} Suppose that $u\in \funcR(\Delta,b)$. Then
 $\det (D^2u)\geq d_1$ everywhere in $\Delta$, where
$$d_1 = \left(\frac{4b
 \mathrm{Diam}(\Delta)^2}{n}\right)^{-n}.$$
\end{lemma}
This
lemma can be found in \cite{D4}.
\begin{lemma}\label{lemma_4.2}
For any $u\in \funcR(\Delta,b),$ we have
$$\|\phi\|_{C^0(M)}
=\|u-v\|_{C^0(\overline{\Delta})},
$$
where $\phi=f-g,$ and $f$ is the Legendre transform of $u$.
\end{lemma}
This fact  was already used by Donaldson in \cite{D6}, and for the proof,
 the readers are referred to the proof of Proposition 4.4 in \cite{SZ}.

\section { Proof of Main Theorem }\label{sect_4}
 \v\n{\bf Proof of Theorem \ref{main_thm}.}
Without loss of generality, we assume that $v$ is normalized at
$0\in\Delta,$ i.e.,
$$v(0)=0,\;\;\nabla v(0)=0. $$
Then $g$ attains its minimum at $0$.
      Suppose that
\begin{equation}\label{eqn_4.2}\|u \|_{C^0(\overline{\Delta})}\leq \mathcal C_1,\;\;\;\|v\|_{C^0(\overline{\Delta})}\leq \mathcal C_1\end{equation}
for some constant $\mathcal C_1>0$. Let   $f$ be the Legendre transforms  of
  $u$. By Lemma \ref{lemma_4.2} and
\eqref{eqn_4.2}, for any $C>4\mathcal C_1$, we have
\begin{equation}\label{eqn_4.3}
S_{g}(C-2\mathcal C_1)\subset S_{f }(C)\subset S_{g}(C+2\mathcal C_1),
\end{equation}
where
$$
S_{f }(C)=\{x|f (x)\leq C\},\;\;\;S_{g}(C)=\{x|g(x)\leq C\}.
$$
\v By   Lemma \ref{lemma_3.2},
\begin{equation}\label{eqn_4.4}
\det (f_{ij})\leq d_1\inv.\end{equation}
On the other hand, by Proposition \ref{prop_3.1} we have
$$
\frac{\det(g_{ij})}{\det (f _{ij})}\leq \mathcal C_2(n,b,\K,\mathcal C_1).
$$
  Restricting to $S_{f }(C )\subset S_{g }(C+2\mathcal C_1 )$, we have
\begin{equation}\label{eqn_4.5}
\det (f _{ij})\geq  {\mathcal C_2^{-1}} \det (g_{ij})\geq d_2(\mathcal C_2,C).
\end{equation}
By the determinant estimates \eqref{eqn_4.4} and  \eqref{eqn_4.5}, applying the Caffarelli-Guti$\acute{e}$rrez theory(see \cite{CG}) and the
Caffarelli-Schauder estimate(see \cite{CC}) in $S_f(C)$, we conclude that
\begin{equation}\label{eqn_4.6}\|f\|_{C^{3,\alpha}(S_g(\frac{C}{2}-2\mathcal C_1))}\leq \|f\|_{C^{3,\alpha}(S_f(\frac{C}{2}))} \leq \mathcal C_3(n, b,\K, \mathcal C_1,C) \end{equation}
for any $C>4\mathcal C_1$(details see \cite{TW}).

\v For any domain $\Omega\subset\subset \Delta,$ by the convexity
of $u $   and \eqref{eqn_4.2} we have
$$\nabla u(\Omega)\subset B_{\frac{2\mathcal C_1}{d}}(0),$$
while the later one is contained in $ S_g(C_ 3)$
for some  constant $C_3>0$.
Here $d=dist(\Omega,\partial \Delta)$ denotes the Euclidean distance from $\Omega$ to $\partial \Delta.$
 Choose $C=2C_3+4\mathcal C_1.$ By \eqref{eqn_4.6},
we obtain $\|u\|_{C^{3,\alpha}(\Omega)} \leq \mathcal C_3(n, b,\K, \mathcal C_1,C) $.   The
Proposition follows from the standard bootstrap argument.
 $\Box$

\end {document}